\magnification\magstep1
\overfullrule = 0pt

\pageno=0
\def\cf{{\it cf. \ }}
\def\ie{{\it i.e. \ }}
\def\ve{\vfill\eject}

\def\n{\noindent}
\def\ms{\medskip}

\def\RR{{\mathop{{\rm I}\kern-.2em{\rm R}}\nolimits}}
\def\TT{{\mathop{{\rm I}\kern-.4em{\rm T}}\nolimits}}
\def\FF{{\mathop{{\rm I}\kern-.2em{\rm F}}\nolimits}}
\def\NN{{\mathop{{\rm I}\kern-.2em{\rm N}}\nolimits}}
\def\EE{{\mathop{{\rm I}\kern-.2em{\rm E}}\nolimits}}
\def\CC{{\rm C\kern-.18cm\vrule width.6pt height 6pt depth-.2pt
\kern.18cm}}
\def\ZZ{{\mathop{{\rm Z}\kern-.28em{\rm Z}}\nolimits}}

\def\sst{\scriptstyle}
\def\qed{{\hfill{\vrule height7pt width7pt
depth0pt}\par\bigskip}} 
\def\vp{\varepsilon}
\def\pf{\medskip\n {\bf Proof.}~~}

\def\nat{\NN}

\centerline{\bf An inequality for p-orthogonal}
\centerline{\bf sums in non-commutative ${\bf L_p}$}
\ms
\centerline{ by}\ms

\centerline {Gilles Pisier\footnote*{Supported in part by the NSF   and
by the Texas Advanced Research Program 010366-163.}}
\centerline{ Texas A\&M University}
\centerline{College Station, TX 77843-3368, USA}
\centerline{and}
\centerline{Universit\'e Paris VI}
\centerline{Equipe d'Analyse, Case 186, 75252}
\centerline{ Paris Cedex 05, France}
\bigskip
\centerline{\it Corrected August 12, 1999.  
}\centerline{\it  Submitted to Illinois J. 
Math.  }
 
\ms \ms \ms 

{\bf Abstract:} We give an alternate proof
of one of the   inequalities proved recently
for martingales (=sums of martingale differences) in a non-commutative
$L_p$-space, with 
  $1<p<\infty$, by Q. Xu and the author. This new approach
is restricted to $p$ an even integer, but it
yields a constant which is $O(p)$ when $p\to \infty$ 
and it applies to a much more general kind of sums
which we call $p$-orthogonal. We use mainly combinatorial tools, namely the
M\"obius inversion formula for the lattice of partitions of a $p$-element set.
  \vfill\eject

\baselineskip = 18pt

\n {\bf \S 0. Introduction}

In a recent paper ([PX]), Quanhua Xu and the author
have proved non-commutative 
versions of the Burkholder-Gundy classical inequalities 
(see [BG, B1-B4]) relating the 
$L_p$-norms of a martingale with those of its square function $(1<p<\infty)$. 
We will continue this investigation here. Our objective is two-fold. First we 
will improve the order of growth of the constant in the main inequality from 
[PX] when $p\to \infty$. We obtain a constant which is $O(p)$ when $p\to 
\infty$, thus yielding the ``sharp'' order of growth. Sharp constants 
themselves are known in the classical-commutative-case, see [B3 and B4, \S
11], but they  seem out of reach of our method.

Secondly, we wish to extend the inequality from martingales to a much broader 
class of sums in non-commutative $L_p$-spaces: \ the $p$-orthogonal sums, 
which are defined as follows.

Let $(M,\tau)$ be a von Neumann algebra equipped with a standard (= faithful, 
normal) trace with $\tau(1)=1$, and let $L_p(\tau)$ be the associated 
``non-commutative'' $L_p$-spaces. (Of course, if $M$ is commutative, we 
recover the classical $L_p$ associated to a probability space.) Let $p\ge 2$ 
be an even integer. A family $d = (d_i)_{i\in I}$ is called $p$-orthogonal if, 
for any {\it  injective} function $g\colon \ [1,2,\ldots, p]\to I$ we have
$$\tau(d^*_{g(1)} d_{g(2)} d^*_{g(3)} d_{g(4)}\ldots d^*_{g(p-1)} d_{g(p)}) = 
0.$$

In the commutative case, \ie for classical random variables, this notion
is very close to that of a ``multiplicative sequence" already considered in
the literature (see Remark 2.4 below for more details).

 Let us assume $I$ finite for simplicity. We
will denote simply by
$\|~~~\|_p$ the norm in
$L_p(\tau)$. 

Then the following inequality, which is our main result, holds:
$$\left\|\sum_{i\in I} d_i\right\|_{L_p(\tau)} \le {3\pi \over 2} p\ S(d,p) 
\leqno (0.1)$$
where we have set
$$S(d,p) = \max\left\{\left\|\left(\sum d^*_id_i\right)^{1/2}\right\|_p, \quad
\left\|\left( \sum d_id^*_i\right)^{1/2}\right\|_p \right\}.\leqno(0.2)$$

Clearly, any 
martingale difference sequence is $p$-orthogonal, but 
the class of $p$-orthogonal sums includes a broader class of sums which appear 
rather naturally in Harmonic Analysis. For instance, let $\Lambda \subset G$ 
be a subset of a discrete group with unit element $e$. We call $\Lambda$ 
$p$-dissociate if for any choice $t_1,t_2,\ldots, t_p$ of $p$ {\it
distinct\/}  points in $\Lambda$ we have
$$t^{-1}_1 t_2 t^{-1}_3 t_4\ldots t^{-1}_{p-1} t_p\ne e.$$
See [Ru] for examples of this in the Abelian case.
Then let $\lambda\colon \ G\to B(\ell_2(G))$ be the left regular 
representation of $G$, let ${\cal M}$ be the von~Neumann algebra generated by 
$\lambda$ and let $\tau_G$ be the usual normalized trace on ${\cal M}$ defined 
by
$$\tau_G(x) = \langle x\delta_e,\delta_e\rangle.$$
Let $(\delta_t)_{t\in G}$ be the canonical basis of $\ell_2(G)$.

With this notation (and with $\tau$ as before), for any function $x\colon \ 
\Lambda \to L_p(\tau)$ the family
$$d_t = \lambda(t) \otimes x(t)$$
is $p$-orthogonal in $L_p(\tau_G\times \tau)$. Therefore (0.1) holds in this 
case too for any finite subset $I\subset \Lambda$. More generally, a family 
$(\Lambda_i)_{i\in I}$ of disjoint subsets of $\Lambda$ will be called 
$p$-dissociate if every family $(t_i)_{i\in I}$ with $t_i\in\Lambda_i$ for all 
$i$ in $I$ is itself $p$-dissociate. Then assuming, say, that $x$ is finitely 
supported, if we define
$$d_i = \sum_{t\in\Lambda_i} \lambda(t) \otimes x(t)$$
we obtain again a $p$-orthogonal sum so that (0.1) holds in this case too. For 
instance in the case $G = \ZZ$ and $\Lambda_i = [2^i,2^{i+1}[$, treating 
separately the cases of $\{\Lambda_i\mid i \hbox{ even}\}$ and $\{\Lambda_i 
\mid i \hbox{ odd}\}$ we can recover from (0.1) one of the classical 
Littlewood-Paley inequalities for Fourier series:
$$\left\|\sum_{n>0} a_ne^{int}\right\|_p \le C_p\|\underline S\|_p$$
where
$$\underline S = \left(\sum_{k\ge 0} \left|\sum_{2^k\le n< 2^{k+1}} a_n 
e^{int}\right|^2 \right)^{1/2}$$
and where, say, we assume that $(a_n)_{n>0}$ is a finitely supported 
sequence of scalars.

A surprising feature of our proof of the martingale inequalities
(or their extensions)
is that we use very elementary tools. Indeed, in the
non-commutative setting which is our main motivation, most of the usual techniques
such as stopping times
or maximal inequalities are unavailable, or apparently inefficient. Therefore,
we must use
only H\"older's inequalities and certain identities.
For example when $p=4$ we are using
an identity of the form:  
 {$$  (\sum d_i)^4=\Sigma
  +6   \sum d_i^4 -8  
(\sum d_i^3)(\sum d_i)
$$
$$- 3   (\sum d_i^2)^2
+ 6 (\sum d_i^2)(\sum d_i)^2 .$$ }
 where $$\Sigma=\sum_{i_1,i_2,i_3,i_4\ 
all\  distinct} 
d_{i_1}d_{i_2}d_{i_3}d_{i_4}   .$$
More generally, for any even integer $p$, there is
an analogous identity
for $(\sum d_i)^p$ in which the coefficients
appearing (such as $6, -8, -3, 6$ when $p=4$) can be explicitly computed using
the M\"obius inversion formula, classical in the combinatorics
of partitions (\cf [R1,R2, A]). In particular, there are explicit
formulae (due to Sch\"utzenberger, see Theorem 1.2 below) for these
coefficients, which lead to suprisingly good bounds for the constants in our
inequalities.

\def\F{{\bf F}}

\n {\bf Remark.0.1.}
 Many examples of non-commutative martingales can be given using
(non-commutative) Harmonic Analysis. Let $G$ be a discrete group, and
let $\lambda_G: \ G\to \ell_2(G)$
be its left regular representation. The von Neumann algebra
of $G$ is defined as $M=\lambda_G (G)''$, and it can be equipped
with the standard trace $\tau$ defined by $\tau_G (x)=\langle x\delta_e, 
\delta_e\rangle $.  Let $G_n$ ($n\in \nat$) be a non-decreasing
sequence of subgroups, and let $M_n=\lambda_G (G_n)''$.
Then, denoting by $\EE_n$ the (contractive) conditional expectation
from $M$ to $M_n$ (which is also contractive
on $L_p(\tau_G)$ whenever $1\le p< \infty$), for
any $f$ in $L_p(\tau_G)$, the sequence $d_n=\EE_nf-\EE_{n-1}f$
is a martingale difference sequence, hence satisfies (0.1).

\n {\bf Remark.0.2.} As explained in [PX], the ``free group filtration" is
a typical example to which the preceding point
  applies. By this we mean the case when     $G=\F_\infty$   the free
group with  countably many generators denoted by $\{ g_0, g_1, g_2,...\}$,  
  $G_n\subset G$ is the subgroup generated
by $\{ g_0,g_1,..,g_n\}$ and again   $M_n=\lambda_G (G_n)''$. We will use
this example below
in one of the proofs. We could consider more generally the filtration
associated to a  free product of a countable collection of groups.

\n {\bf Remark.0.3.} Another example is the ``free-Gaussian" analog of the
preceding. Let
$(M,\tau)$ be
a von Neumann algebra equipped with a standard  normalized trace. Let
$(x_n)_{n\ge 0}$ be a free
semi-circular family in Voiculescu's sense ([VDN]) in $(M,\tau)$, and
let $M_n$ be the von Neumann algebra generated
in $M$ by $\{x_0,...,x_n\}$.
Then again $(M_n)$ is an interesting example to which 
(0.1) applies. This case was recently studied by Biane and Speicher
 [BS].
 Their main result gives evidence that,
 for martingales relative to the free group filtration and its free-Gaussian
analog, 
  the constant appearing
 in  (0.1) might actually be
 {\it bounded} when $p\to \infty$, but this remains open.

 \n{\bf Acknowledgement.} I thank the referee for his careful reading of the
manuscript.
\ve

\n {\bf \S 1. M\"obius Inversion}

We will use crucially some well known ideas from the combinatorial theory of 
partitions, which can be found, for instance, in Rota's texts ([R1, R2]) or in 
the book [A]. We denote by $P_n$ the lattice of all partitions of $[1,\ldots, 
n]$, equipped with the following order:\ we write $\sigma\le \pi$ (or 
equivalently $\pi\ge\sigma$) when every ``block'' of the partition $\sigma$ is 
contained in some block of $\pi$. Let $\dot 0$ and $\dot 1$ be respectively 
the minimal and maximal elements in $P_n$, so that $\dot 0$ is the partition 
into $n$ singletons and $\dot 1$ the partition formed of the single set 
$\{1,\ldots, n\}$. We denote by $\nu(\pi)$ the number of blocks of $\pi$ (so
that $\nu(\dot 0) = n$ and $\nu(\dot 1) = 1$).

For any $\pi$ in $P_n$ and any $i=1,2,\ldots, n$, we denote by $r_i(\pi)$ the 
number of blocks (possibly $=0$) of $\pi$ of cardinality $i$. In particular, 
we have $\sum^n_1 ir_i(\pi) = n$ and $\sum^n_1 r_i(\pi) = \nu(\pi)$.

Given two partitions $\sigma,\pi$ in $P_n$ with $\sigma\le \pi$ we denote by 
$\mu(\sigma,\pi)$ the M\"obius function, which has the following fundamental 
property.

\proclaim Proposition 1.1. Let $V$ be a vector space. Consider two functions 
$\Phi\colon \ P_n\to V$ and $\Psi\colon \ P_n\to V$.
\item{(i)} If 
$$\eqalignno{\Psi(\sigma) &= \sum\limits_{\pi\le \sigma} \Phi(\pi),\cr
\noalign{\hbox{then}}
\Phi(\sigma) &= \sum_{\pi\le\sigma} \mu(\pi,\sigma)\Psi(\pi).}$$
\item{(ii)} If
$$\eqalignno{\Psi(\sigma) &= \sum_{\pi\ge\sigma} \Phi(\pi),\cr
\noalign{\hbox{then}}
\Phi(\sigma) &= \sum_{\pi\ge \sigma} \mu(\sigma,\pi) \Psi(\pi).}$$
\item{(iii)} In particular we have:
$$\sum_{0\le \pi \le \sigma} \mu(\pi,\sigma) = 0.\leqno \forall~\sigma\ne \dot 
0$$

\n {\bf Remark.} (iii) follows from (i) applied with $\Phi$ equal to the delta 
function at $\dot 0$ (i.e.\ $\Phi(\pi) = 0$ $\forall~\pi\ne \dot 0$ and 
$\Phi(\dot 0) = 1$) and $\Psi\equiv 1$.

We also recall Sch\"utzenberger's theorem (see [A] or [R1-2]):

\proclaim Theorem 1.2. For any $\pi$ we have
$$\mu(\dot 0, \pi) = \prod^n_{i=1} [(-1)^{i-1} (i-1)!]^{r_i(\pi)},$$
and consequently
$$\sum_{\pi\in P_n} |\mu(\dot 0,\pi)| = n!.\leqno (1.1)$$

We now apply these results to set the stage for the questions of interest to 
us. Let $E_1,\ldots, E_n, V$ be vector spaces equipped with a multilinear
form  (= a ``product'')  
$$\varphi\colon \ E_1\times\cdots\times E_n\to V.$$
 Let $I$ 
be a finite set. For each $k=1,2,\ldots, n$ and $i\in I$, we give ourselves 
elements $d_i(k) \in E_k$, and we form the sum
$$F_k = \sum_{i\in I} d_i(k).$$
Then we are interested in ``computing'' or ``expanding'' in a specific manner 
the quantity
$$\varphi(F_1,\ldots, F_n).$$
We can start by writing obviously
$$\varphi(F_1,\ldots, F_n) = \sum_g \varphi(d_{g(1)}(1),\ldots, d_{g(n)}(n))$$
where the sum runs over all functions $g\colon \ [1,2,\ldots, n] \to I$. Let 
$\pi(g)$ be the partition associated to $g$, namely the partition obtained 
from $\bigcup\limits_{i\in I} g^{-1}(\{i\})$ after deletion of all the empty 
blocks.
We can write
$$\varphi(F_1,\ldots, F_n) = \sum_{\sigma\in P_n} \Phi(\sigma)$$
where
$\Phi(\sigma) = \sum_{g\colon\ \pi(g) =\sigma} \varphi(d_{g(1)}(1),\ldots, 
d_{g(n)}(n))$.

\n 
By Theorem 1.1, if we let $\Psi(\sigma) = \sum\limits_{\pi\ge\sigma} 
\Phi(\pi)$ we can write using (ii) and (iii) in Proposition~1.1:
$$\eqalign{\varphi(F_1,\ldots, F_n) &= \Phi(\dot 0) + \sum_{\dot 0<\sigma} 
\Phi(\sigma)\cr
&= \Phi(\dot 0) + \sum_{\dot 0<\sigma} \sum_{\pi\ge \sigma} \mu(\sigma,\pi) 
\Psi(\pi)\cr
&= \Phi(\dot 0) + \sum_{\dot 0<\pi} \Psi(\pi)\cdot \sum_{\dot 0<\sigma\le\pi} 
\mu(\sigma,\pi)\cr
&= \Phi(\dot 0) - \sum_{\dot 0<\pi} \Psi(\pi) \mu(\dot 0,\pi).}$$
Recapitulating, we state:

\proclaim Corollary 1.3. The following identity holds
$$\varphi(F_1,\ldots, F_n) = \Phi(\dot 0) - \sum_{0<\pi} \Psi(\pi) \mu(\dot 0, 
\pi)$$
where
$$\Phi(\dot 0) = \sum_{g~{\rm injective}} \varphi(d_{g(1)}(1),\ldots, 
d_{g(n)}(n))$$
and where
$$\Psi(\pi) = \sum_{g\colon \ \pi(g)\ge \pi} \varphi(d_{g(1)}(1),\ldots, 
d_{g(n)}(n)).$$
\ve

\n {\bf \S 2. The Commutative Case}

Although the main point of this paper is the non-commutative case, we prefer to 
present the proof first in the classical setting. This will make it much 
easier for the reader to follow the arguments in the next sections. Note that 
although many results similar to our Theorem~2.1 below exist in the literature 
(\cf e.g. [St] and Remark 2.4 below), we could not quite find a reference
for the same result.

Let $(\Omega,m)$ be any measure space and let $p=2k$ be an even integer. Let 
$(d_i)_{i\in I}$ be a finite sequence in $L_p = L_p(m)$. We will say that 
$(d_i)_{i\in I}$ is $p$-{\it orthogonal\/} if for any injective map $g\colon \ 
I\to [1,\ldots, p]$ we have
$$\int \bar d_{g(1)} d_{g(2)} \bar d_{g(3)}\ldots \bar d_{g(p-1)} d_{g(p)} dm 
=0.\leqno (2.1)$$
Clearly, if $p=2$ we recover the usual orthogonality in $L_2$. Throughout this 
section, we will denote
$$S = \left(\sum_{i\in I} |d_i|^2\right)^{1/2}.$$
It is easy to check that any martingale difference sequence in $L_p$ is 
$p$-orthogonal (consider the largest value of $g$, say $g(i) = n$ and take the 
conditional expectation of index $n-1$, before the integral in (2.1)).

\proclaim Theorem 2.1. Let $(d_i)_{i\in I}$ be a $p$-orthogonal finite 
sequence in $L_p = L_p(\Omega, m)$. We have then for all even integers 
$p=2k$:
$$A_p\|S\|_p - \left(\sum_{i\in I} \|d_i\|^p_p\right)^{1/p} \le 
\left\|\sum_{i\in I} d_i\right\|_p \le 2p \|S\|_p\leqno (2.2)$$
where $0<A_p\le 1$ is a constant depending only on $p$.

It is well known that a random variable $f$ on a probability space is
exponentially integrable,
\ie 
$$\exists\  \delta >0  \ \hbox{such that}\  \int \exp(\delta |f|) dP<\infty$$
iff
$f\in L_p$ for any even   integer $p>0$ and 
$$\exists\  K \ \hbox{such that}\   \|f\|_p
\le K p   \quad \forall p>0  \ \hbox{even integer}.
$$ Moreover, the corresponding norms are equivalent. Thus we have:

\proclaim Corollary 2.2. Let $(d_i)_{i\in I}$ be a (countable)
family of random variables on a probability space $(\Omega,P)$  which are
$p$-orthogonal
  for any even integer $p=2k$. Then,  
if the ``square function"
$S=(\sum |d_i|^2)^{1/2}$ is in the unit ball of  $L_\infty$,
we have
$$   \int \exp(\delta |\sum d_i |) dP\le 2 $$
where $\delta>0$ is a numerical constant (independent of the family $(d_i)$).
 
\n{\bf Proof of Theorem 2.1.} For simplicity we restrict ourselves to the
$\RR$-valued case. We apply  the combinatorics in \S 1 to the multilinear
form:
$$\varphi\colon \ L_p \times   \cdots \times L_p \to \RR$$
defined by $\varphi(x_1,x_2,\ldots, x_{p-1},x_p) = \int x_1x_2\ldots x_{p-1} 
x_p\ dm$. The hypothesis in Theorem~2.1 guarantees that $\Phi(\dot 0) = 0$. Let 
$f = \sum\limits_{i\in I} d_i$. Applying Corollary 1.3, we thus obtain:
$$\|f\|^p_p = - \sum_{\dot 0<\pi} \mu(\dot 0,\pi) \Psi(\pi)\leqno (2.3)$$
where
$$\Psi(\pi) = \int \prod^p_{j=1} \left(\sum_{i\in I} d^j_i\right)^{r_j(\pi)} 
d\mu.$$
If $ j\ge 2$, then  $\left|\sum d^j_i\right|^{1/j} \le S$,  so by 
H\"older's inequality $\left({r_1(\pi)\over p} + {p-r_1(\pi)\over p}=1\right)$
$$ {|\Psi(\pi)| \le \int |f|^{r_1(\pi)} S^{p-r_1(\pi)}d\mu 
\le \|f\|^{r_1(\pi)}_p \|S\|^{p-r_1(\pi)}_p.}$$
Thus we obtain
$$\|f\|^p_p \le \sum_{\dot 0 <\pi} |\mu(\dot 0,\pi)|\ \|f\|^{r_1(\pi)}_p 
\|S\|^{p-r_1(\pi)}_p.$$
Note that ${\dot 0 <\pi}$ implies  ${r_1(\pi)}\le p-2$, hence the last sum can
be rewritten as
$$\sum_{0\le r\le  p-2} \|f\|^r_p \|S\|^{p-r}_p a_r \hbox{ with } a_r = 
\sum_{r_1(\pi)=r} |\mu(\dot 0,\pi)|.$$
A moment of thought shows that $a_r= {p\choose r} b_r$ where $b_r$ is the sum 
of $|\mu(\dot 0, \sigma)|$ over all partitions $\sigma$ of $[1,\ldots, p-r]$ 
without any singleton. A fortiori by (1.1), we have 
$$b_r\le (p-r)!$$
Thus we obtain finally
$$\|f\|^p_p \le \sum_{0\le r\le  p-2} \|f\|^r_p \|S\|^{p-r}_p {p\choose r}
(p-r)!$$ Therefore, using the sublemma below, we conclude that
$$\|f\|_p \le 2p\|S\|_p.$$

\proclaim  Sublemma 2.3. Let $x,y$ be positive numbers such that
$$x^p \le \sum_{0\le r<p} x^r y^{p-r} {p\choose r} (p-r)!$$
Then $x\le 2py$.

\pf Let $t=y/x$. We have
$$1 \le \sum_{0\le r<p} {p\choose r} t^{p-r} (p-r)!.$$
Using $\int^\infty_0 s^{p-r} e^{-s} ds = (p-r)$! and $\int^\infty_0 e^{-s} ds 
= 1$, we obtain
$$1 \le \int^\infty_0 [(1+ts)^p -1] e^{-s}ds = \int^\infty_0 (1+ts)^p e^{-s} 
ds - 1$$
whence $2 \le  \int^\infty_0 \exp(pts-s)ds$.

Therefore if $pt<1$ this implies $2\le (1-pt)^{-1}$ hence ${1\over t} \le 2p$ 
(and if $pt\ge 1$, then ${1\over t} \le p$ which is even better).\qed

We now turn to the converse inequality.

With the same notation as before, we now ``isolate'' in (2.3) the terms 
corresponding to the partitions $\pi$ such that $r_2(\pi) = p/2$, i.e.\ $\pi$ 
is a partition of $[1,\ldots, p]$ into $p/2$ pairs. Let $\alpha_p$ be the 
number of such partitions.
For such a $\pi$, by Theorem 1.2 we have $\mu(\dot 0,\pi) = (-1)^{p/2}$ and 
$\Psi(\pi) = \int \left(\sum|d_i|^2\right)^{p/2}dm = \|S\|^p_p$. Thus we 
obtain
$$\|f\|^p_p = \alpha_p(-1)^{p/2+1} \|S\|^p_p - \sum{}' \mu(\dot 0, \pi) 
\Psi(\pi) \leqno (2.4)$$
where the symbol $\sum'$ means that we sum over all $\pi$ with $r_1(\pi) \le 
p-2$  and $r_2(\pi) < p/2$. A simple calculation shows that
$$\alpha_p = p![2^{p/2}(p/2)!]^{-1}.$$
We can write
$$\sum{}' \mu(\dot 0,\pi) \Psi(\pi) = \sum_{0\le r\le  p-2} C(r) \leqno
(2.5)$$ where $C(r) = \sum\limits_{\sst r_1(\pi) = r\atop\sst r_2(\pi)<p/2}
\mu(\dot  0,\pi) 
\Psi(\pi)$.

By arguing as above, we obtain
$$|C(r)| \le {p\choose r} (p-r)! \|f\|^r_p S^{p-r}.$$
 But now, this estimation 
will be sufficiently efficient for our purposes only if $r>0$; the term 
$C(0)$ has to be estimated separately. We have
$$\eqalign{|C(0)| &\le \sum_\lambda \hbox{card}(\pi\mid r_j(\pi) =
\lambda_j, \forall j\ge 0) 
\Pi((i-1)!)^{\lambda_i}\cr
&\quad \cdot \int\left(\sum d^2_i\right)^{\lambda_2} \left(\sum 
|d_i|^3\right)^{\lambda_3} \ldots \left(\sum |d_i|^p\right)^{\lambda_p} dm,}$$
where the sum runs over all integers $\lambda_j\ge 0$ such that $p = \lambda_1 
+ 2\lambda_2 + \cdots+ p\lambda_p$ with $\lambda_2 < p/2$ and $\lambda_1=0$.

Since $2<3\le p$, we can write ${1\over 3} = {1-\theta\over 2} + {\theta\over
p}$  with $\theta>0$. Hence, by H\"older:
$$\left\|\left(\sum |d_i|^3\right)^{1/3}\right\|_p \le \|S\|^{1-\theta}_p 
\left(\sum \|d_i\|^p_p\right)^{\theta/p}.$$
Let $h = \left(\sum\limits_{i\in I} \|d_i\|^p_p\right)^{1/p}$. Since 
$\lambda_2 <p/2$, we have $2\lambda_2 \le p-2$ and since we may as well
assume $h \le \|S\|_p$ (otherwise the left side of (2.2) is negative), we
obtain  again by H\"older:
$$\eqalign{\int\left(\sum |d_i|^2\right)^{\lambda_2} \ldots \left(\sum 
|d_i|^p\right)^{\lambda_p} dm &\le \|S\|^{2\lambda_2}_p \left\|\left( 
\sum|d_i|^3 \right)^{1/3}\right\|^{p-2\lambda_2}_p\cr
&\le \|S\|^{2\lambda_2 +(1-\theta)(p-2\lambda_2)}_p \cdot 
h^{\theta(p-2\lambda_2)}  \le \|S\|^{p-2\theta}_p h^{2\theta}.}$$
Thus, returning to (2.4) and (2.5) we can write
$$\|f\|^p_p \ge \alpha_p\|S\|^p_p - \sum_{0<r<p-1} |C(r)| - |C(0)|$$
which implies
$$\alpha_p\|S\|^p_p \le \|f\|^p_p + \sum_{0<r<p-1} {p!\over r!} 
\|f\|^r_p S^{p-r}
 + \beta_p \|S\|^{p-2\theta}_p h^{2\theta}$$
where $\beta_p$ is a constant depending only on $p$. Clearly, since
$\theta>0$,  this last  estimate shows that the ratio $\|S\|_p \cdot
[\max\{\|f\|_p, h\}]^{-1}$ must  be bounded above by a constant depending only
on $p$. This yields the left  side of (2.1). \qed
\ms

\n {\bf Remark 2.4.} The literature contains numerous attempts
to generalize orthogonality. For instance, in Stout's book [St] a sequence
of (real valued) random variables is called ``multiplicative"
(resp. ``multiplicative of order $r$")
if it admits moments of all order (resp. of all order $\le r$) and is
$p$-orthogonal for all $p$ (resp. for all $p \le r$).
We are aware of  works by Azuma (1967), Serfling (1969), 
Dharmadhikari and Jogdeo (1969) (for which we refer
to  [St] for precise references) which all relate to
the notion of $p$-orthogonality, but we could not
find results like Theorem 2.1 in the literature, although it
might be known.  One notable exception is the paper [Se] 
(see also [LS]) which contains a statement ([Se, Th. 2.1]) similar to  the
right side of (2.2), namely it is proved there that there is a constant $A$
such that for any $p$-orthogonal family
$(d_i)_{i\in \nat}$ and any $n$, we have
$$\|\sum_1^n d_i\|_p\le A n^{1/2} \sup_{i\in \nat} \|d_i\|_p .\leqno(2.6)$$
 Note that (2.6)  follows also from the right side of (2.2). The 
 basic idea of
 the proof of (2.6) in [Se] turns out to be
 essentially the same as the one used above
 for the right side of
(2.2), but the dependence of
 $A$ with respect to $p$ (or the connection
 with the combinatorics of partitions) does not appear in [Se].
 (I am very grateful to Prof. Serfling for kindly communicating
 to me a copy of this paper upon request, to allow
 a comparison with the above results.)

\n {\bf Remark 2.5.} As a corollary, we obtain a proof of the classical 
Burkholder-Gundy inequalities, which say that $\|S\|_p$ and $\left\|\sum 
d_n\right\|_p;$ are equivalent whenever $d = (d_n)$ is a martingale difference 
sequence. Indeed, as already mentioned, these are $p$-orthogonal. Moreover, 
the inequality $\left(\sum \|d_n\|^p_p\right)^{1/p} \le 2\left\|\sum d_n 
\right\|_p$ is elementary (by interpolation between $p=2$ and $p=\infty$). 
Therefore, (2.2) implies in this case that for any choices of signs $\vp_n = 
\pm 1$ we have
$$\left\|\sum \vp_nd_n\right\|_p \le C_p\left\|\sum d_n\right\|_p.\leqno 
(2.7)$$
Finally interpolation and duality starting from (2.7) allow to pass from $p$ 
an even integer to the whole range $1<p<\infty$.

\n Note there is a well known very classical proof 
 due to Paley [Pa] (for dyadic martingales),
which also is based on the case when $p$ is an even integer, but Paley's
proof uses the ``martingale assumption" several times
(and not merely $p$-orthogonality), moreover
he uses the maximal inequalities, which do not seem
to have a counterpart for non-commutative martingales.\ms

\n {\bf Remark 2.6.} Note that we cannot have a lower bound $A_p\|S\|_p \le 
\left\|\sum d_i\right\|_p$ for general $p$-orthogonal sums. Indeed, just 
taking a pair $d_1,d_2$ and the rest equal to zero, we see that this would 
imply when $p=4$ that $\|d_1\|_p \le A^{-1}_p \|d_1+d_2\|_p$ which is clearly 
absurd without any assumption of the pair $d_{1},d_2$. (Note in particular that
$p$-orthogonality does not even imply linear independence!)

\ve
\n {\bf \S 3. The non-commutative Case}

Let $M$ be a von Neumann algebra equipped with a faithful normal and
normalized
trace
$\tau$.  Let $1\le p <\infty$. The space $L_p(M,\tau)$ (or simply
$L_p(\tau)$) is  defined as the completion of $M$ with respect to the norm
$\|x\|_p = 
\tau(|x|^p)^{1/p}$ (here of course $|x| = (x^*x)^{1/2}$). It is natural to 
set, say by convention, $L_\infty(\tau) = M$ equipped with the operator norm.

Now if $p$ is an even integer we say that a finite sequence $(d_i)_{i\in I}$ in 
$L_p(\tau)$ is $p$-orthogonal if, for any {\it injective\/} map $g\colon \ 
[1,\ldots, p] \to I$, we have
$$\tau(d^*_{g(1)} d_{g(2)}\ldots d^*_{g(p-1)} d_{g(p)}) = 0.$$
Observe that $p$-orthogonality is inherited by subfamilies,
and also, 
that if the cardinality of $I$ is
$<p$ any family $d = (d_i)_{i\in I}$ is   $p$-orthogonal,
but this is actually irrelevant for our
purposes, since we are only interested in the case
when $I$ is large compared with $p$.

Of course if $M$ is commutative, then $(M,\tau)$ can be identified with 
$L_\infty(\Omega,m)$ for some measure space $(\Omega,m)$ and $\tau(x) =  \int 
x dm$, so that we recover the notion introduced in the preceding 
section. The main result of this paper is the following non-commutative 
version of Theorem~2.1.

\proclaim Theorem 3.1. Let $(M,\tau)$ be as above. Let $p>2$ be an even 
integer. Then for any $p$-orthogonal finite sequence $(d_i)_{i\in I}$ in 
$L_p{(\tau)}$, we have
$$ \left\|\sum_{i\in I} d_i\right\|_{L_p(\tau)} \le 
{3\pi\over 2}p\|S\|_{L_p} \leqno (3.1)$$
 where the ``square 
function'' $S$ is defined as
$$S = \left(\sum_{i\in I} d^*_id_i + d_id^*_i\right)^{1/2}. \leqno (3.2)$$
In particular, when $I$ is infinite, if  $S$ converges
(strong operator topology) to a bounded operator in the unit ball of  $M$,
and if $(d_i)_{i\in I}$ is $p$-orthogonal for all $p$,
then the series $\sum d_i$ obviously converges in ${L_2(\tau)}$ and
its sum satisfies
$$\tau\left(\exp(\delta |\sum d_i|)\right) \le 2$$
where $\delta>0$ is a numerical constant (independent of the family $(d_i)$).

\pf Let $f = \sum\limits_{i\in I} d_i$. We can write as before
$$\tau[(f^*f)^{p/2}] = -\sum_{\dot 0<\pi} \mu(\dot 0,\pi) \Psi(\pi)$$
where $\Phi$ and $\Psi$ are now defined as follows:
$$\eqalign{\Phi(\sigma) &= \sum_{g\colon \ \pi(g)=\sigma} \tau(d^*_{g(1)} 
d_{g(2)} \ldots d^*_{g(p-1)} d_{g(p)})\cr
\Psi(\pi) &= \sum_{\sigma\ge \pi} \Phi(\sigma) \hbox{ or equivalently}\cr
\Psi(\pi) &= \sum_{g\colon \ \pi(g)\ge \sigma} \tau(d^*_{g(1)} d_{g(2)} \ldots 
d^*_{g(p-1)} d_{g(p)}).}$$
A quick inspection of the proof of Theorem~2.1 shows that all we need is the 
next statement.\ms

\proclaim   Sublemma 3.2. For any partition $\pi$, we have
$$|\Psi(\pi)| \le (\alpha\|S\|_p)^{p-r_1(\pi)} \|f\|^{r_1(\pi)}_p,$$
where $\alpha= {3\pi\over 4}$.

\n Indeed, using this and arguing as for Theorem~2.1, we obtain
$$\|f\|^p_p \le \sum_{0\le r<p} \|f\|^r_p (\alpha\|S\|_p)^{p-r} {p\choose r} 
(p-r)!,$$
hence by Sublemma 2.3, we conclude that
$$\|f\|_p \le 2\alpha p\|S\|_p.\eqno $$
This shows (3.1). The last assertion in Theorem 3.1 is then deduced
from this
exactly like Corollary 2.2 was deduced from Theorem 2.1. We leave the details
to the reader.\qed

\n {\bf Remark.} The inequality (3.1) probably admits a converse (analogous
to the left side
of (2.1)), but we could not prove it.
The difficulty lies in the fact that (when, say, 
$p=4$) terms such as 
$$\psi = \sum_{ij} \tau(d^*_id_jd^*_id_j)$$
may be negative in the non-commutative case. For instance, if $(d_i)_{1\le i\le
n}$ is a family of anti-commuting self-adjoint unitaries (= a spin system) then
$d^*_i d_jd^*_id_j = -I$ for all $i\ne j$ and it is equal to $I$ otherwise.
Hence, in this case $\psi = n-(n^2-n) = 2n-n^2$.

\n {\bf Remark.} The above proof actually shows that
$\|f\|_p \le 2\alpha p S(d,p),$
with $S(d,p)$ as defined in (0.2).

To prove Sublemma 3.2, we need several more lemmas. In the first one, we
denote by $\F_I$ the free group with free generators
$(g_i)_{i\in I}$ and by $\varphi$ the normalized  trace on the von Neumann
algebra of $\F_I$ (essentially as in Remark 0.2). \ms

\proclaim  Sublemma 3.3.  Fix $p\ge 2$ and let $\pi \in P_p$. Let $B_1$ be
the union
of all the singletons of $\pi$, and let $B_2$ be the complement 
of $B_1$ in $[1,\ldots, p]$. Let
$f_k = 
\sum\limits_{i\in I} d_i(k)$ be a (finite) sum in $L_p(\tau)$. Let $\tilde f_k
= 
\sum\limits_{i\in I} \lambda(g_i) \otimes d_i(k)$ in $L_p(\varphi \times 
\tau)$. Then, for a suitable discrete group $G$, there are elements 
$F_1,\ldots, F_p$ in $L_p(\tau_G\times \tau)$ satisfying
$$\forall k \in B_2 \quad \|F_k\|_p = \|\tilde f_k\|_p\quad \hbox{  and 
}\quad
 \forall k\in B_1 \quad \|F_k\|_p = \|  f_k\|_p, \leqno (3.3)$$   and such
that
$$\sum_{\pi(g)\ge \pi} \tau (d_{g(1)}(1) \ldots d_{g(p)}(p)) = (\tau_G
\otimes 
\tau )[F_1F_2\ldots F_p]. \leqno (3.4)$$

\pf Consider first the case when $\pi$ has only one block $[1,\ldots, p]$, 
i.e.\ we want to rewrite
$$\psi = \sum_{i\in I} \tau(d_i(1)\ldots d_i(p)).$$
Then if $p=2$ this is easy, we can write
$$\eqalign{\psi &= \sum_{i,j\in I} \varphi(\lambda(g_i)^* \lambda(g_j)) 
\tau(d_i(1) d_j(2))\cr
&= (\varphi\times\tau)[F_1F_2] \hbox{ where}\cr
F_1 &= \sum \lambda(g_i)^* \otimes d_i(1),\quad F_2 = \sum \lambda(g_j) 
\otimes d_j(2),}$$
and we obtain the announced result. 

\n Assume now that $\pi$ has one block $[1,\ldots, p]$ but that $p$ is 
arbitrary. Let
$$\eqalign{F_1 &= \sum \lambda(g_i)^* \otimes 1\otimes 1\ldots 1 \otimes 
d_i(1)\cr
F_2 &= \sum \lambda(g_i) \otimes \lambda(g_i)^* \otimes 1\cdots \otimes 1 
\otimes d_i(2)\cr
F_3 &= \sum 1 \otimes \lambda(g_i) \otimes \lambda(g_i)^* \otimes\cdots\otimes 
1 \otimes d_i(3)}$$
and so on, until
$$\eqalign{F_{p-1} &= \sum 1 \otimes \cdots \otimes 1 \otimes \lambda (g_i) 
\otimes \lambda(g_i)^* \otimes d_i(p-1)\cr
F_p &= \sum 1\otimes\cdots \otimes 1 \otimes \lambda(g_i) \otimes d_i(p).}$$
Then it is easy to check that (3.3) holds in $L_p(\tau_G\times\tau)$ where $G$ 
is a product of suitably many copies of the free group $\FF_I$. Moreover, we 
clearly have $\psi = (\tau_G \otimes\tau) [F_1F_2\ldots F_p]$. In addition, we 
have produced a group $G$ and families $(\xi^1_i)_{i\in I},\ldots, 
(\xi^p_i)_{i\in I}$ in $VN(G)$ such that, for any map $g\colon \ [1,\ldots, p] 
\to I$, we have  $\tau_G(\xi^1_{g(1)}\ldots \xi^p_{g(p)}) \ne 0$ if and only
if
$g(i) =  g(j)$ $\forall~i,j$ and in that case the non-zero value is equal to
1.

It is now easy to see the recipe for the general case.

Let $A_1,\ldots, A_\nu$ be the blocks of the partition $\pi$ with more than 
one element. We will introduce discrete groups $G_1,\ldots, G_\nu$ and their 
product $G = G_1\times\cdots\times G_\nu$. Let $VN(G)$ denote the von~Neumann 
algebra of $G$, generated by the left regular representation $\lambda_G$. We 
will identify $\lambda_G$ with $\lambda_{G_1}\otimes\cdots\otimes 
\lambda_{G_\nu}$ and $VN(G)$ with $VN(G_1) \overline\otimes \cdots 
\overline\otimes VN(G_\nu)$. 

\n For each $q$ with $1\le q\le \nu$ the previous argument (applied to each
block  separately) produces elements $(\xi^q_i)_{i\in I}$ in $VN(G_q)$ such
that for  any function $g\colon \ A_q\to I$,
$\tau_{G_q}\left(\prod\limits_{a\in A_q} 
\xi^q_{g(a)}\right) = 1$ iff $g$ takes one single value only and $=0$ 
otherwise. (Here the product sign is meant to respect the order of the 
elements in $A_q$.)

\n Then we define
$$F_k\in VN(G_1) \otimes\cdots\otimes VN(G_\nu) \otimes L_p(\tau)$$
as follows:
$$\leqalignno{F_k &= \sum_{i\in I} \xi^1_i \otimes 1 \otimes \cdots\otimes 1 
\otimes d_i(k)& \forall~k\in A_1\cr
F_k &= \sum_{i\in I} 1 \otimes \xi^2_i \otimes 1\cdots \otimes d_i(k) 
&\forall~k\in A_2\cr
F_k &= \sum_{i\in I} 1 \otimes \cdots\otimes \xi^\nu_i \otimes d_i(k). 
&\forall~k\in A_\nu}$$
Finally, if $k\notin A_1\cup\cdots \cup A_\nu$ (i.e.\ $k$ belongs to some 
singleton block of the partition $\pi$) we set 
$$F_k = 1\otimes\cdots \otimes 1 \otimes f_k.$$
It is then easy to check that (3.4) holds. Finally going back to the definition 
of $(\xi^q_i)_{i\in I}$ we see that (3.3) holds. Indeed, it is well known that
we have (this is analogous to Fell's absorption principle)
$$\left\|\sum \lambda(g_i)^* \otimes \lambda(g_i) \otimes 1 \otimes d_i\right\|_p = 
\left\|\sum \lambda(g_i) \otimes d_i\right\|_p.$$
The latter identity can be checked easily in our case by expanding
the $p$-th powers of the sums on
both sides and observing that the corresponding moments are 
pairwise identical. We leave this to the reader.
\qed
\proclaim Lemma 3.4. For any $p\ge 2$ even integer, we have for any $d = 
(d_i)_{i\in I}$ in $L_p(\tau)$
$$\left\|\sum_{i\in I} \lambda(g_i) \otimes  d_i \right\|_p \le {3\pi\over 4}
S(d,p).\leqno  (3.5)$$

We will deduce this from the next result.
The inequality (3.6) below is   due to Buchholz [Bu2], we include
a slightly different argument (and (3.7) is well known).

\proclaim Lemma 3.5. Let $p\ge 2$ be an even integer. Let $(c_i)_{i\in I}$ be a 
free circular family in Voiculescu's sense (\cf  [VDN]) normalized so that 
$\varphi(|c_i|^2)=1$ and $\|c_i\|_\infty=2$. Then we have for all $d = 
(d_i)_{i\in I}$ in $L_p(\tau)$
$$\left\|\sum_{i\in I} c_i\otimes d_i\right\|_p \le K_p S(d,p)\leqno (3.6)$$
where $K_p = \left[{p\choose p/2} {1\over 1+p/2}\right]^{1/p}\le 2$.
Moreover, we also have
$$\left\|\sum_{i\in I} \lambda(g_i) \otimes  d_i\right\|_p \le
3\pi/8\left\|\sum  c_i\otimes d_i\right\|_p.\leqno (3.7)$$

\pf Let $p=2q$. By [Sp1] (see also [BSp] and [HT]), we know that we can write
$$\left\|\sum  c_i\otimes d_i\right\|^p_p = \sum_{\pi\in S^{nc}_q}
\sum_{i_1i_2 
\ldots i_q\in I} \tau(d^*_{i_1} d_{i_{\pi(1)}} \ldots d^*_{i_q} 
d_{i_{\pi(q)}})$$
where the first sum runs over a certain subset $S^{nc}_q$ of the set $S_q$ of 
all permutations of $[1,\ldots, q]$. This subset is defined as follows. We 
consider the sequence of numbers $\Omega=  [1, \pi(1), 2, \pi(2), \ldots, 
q,\pi(q)]$. We will associate to $\pi$ a partition of $[1,2,\ldots, 2q]$ into 
disjoint pairs like this:
Let $1\le i < j\le 2q$. Then we say that the two-point set $[i,j]$ belongs to
the partition
 if, in $\Omega$, we find the same  number at both the $i$-th and the
$j$-th place.  Clearly this is indeed a partition of $[1,\ldots, 2q]$ into
pairs composed of an  odd and an even integer. We will denote by $S^{nc}_q$
the set of permutations 
$\pi$ such that the associated partition just defined is non-crossing (\cf
[K, Sp2]). It can be  shown by a counting argument (\cf  [K]) that
$\hbox{card}(S^{nc}_q) = {2q\choose  q} {1\over q+1}$ (Catalan number). Hence
we have
$$\left\|\sum  c_i\otimes d_i\right\|_p \le K_p \gamma$$
where $\gamma$ is the positive number defined by
$$\gamma^p = \max_{\pi \in S^{nc}_q} \left\{\left| \sum_{i_1,\ldots, i_q} 
\tau(d^*_{i_1} d_{i_{\pi(1)}} \ldots d^*_{i_q}
d_{i_{\pi(q)}})\right|\right\}.$$ Thus the proof of (3.6) can be easily
completed using Lemma~3.6 below (perhaps of some independent interest). To
check  (3.7) we can note that by Voiculescu's results, the family
$(c_i)_{i\in I}$ has  the same distribution as a family of the form
$(u_i|c_i|)_{i\in I}$ where 
$(u_i)_{i\in I}$ and $(c_i)_{i\in I}$ are $*$-free  and where
$(u_i)_{i\in I}$  and $(\lambda(g_i))_{i\in I}$ have the same
$*$-distribution (in the sense
of [VDN]). Let
$\delta = 
\varphi(|c_i|)$ (independent of $I$). A simple computation shows that
$ \delta =8/3\pi$. In addition, note 
that $(u_i\varphi(|c_i|))_{i\in I}$ can be viewed as obtained by a suitable 
conditional expectation from $(u_i|c_i|)_{i\in I}$. Hence we can write:
$$\delta\left\|\sum u_i\otimes d_i\right\|_p \le  \left\|\sum u_i|c_i|
\otimes
d_i
\right\|_p = \left\|\sum c_i\otimes d_i\right\|_p$$
which yields (3.7). \qed

\n{\bf Proof of Sublemma 3.2.} We apply H\"older's inequality to the right
side of (3.4), then we use (3.3) and (3.5) to obtain Sublemma 3.2. \qed

\proclaim Lemma 3.6. Let $(d_i(k))_{i\in I}$, $k=1,2,\ldots, p$ (with $p=2q$ as 
above) be families of elements in $L_p(\tau)$. Then, for all $\pi$ in 
$S^{nc}_p$, we have
$$\left|\sum_{i_1i_2\ldots i_q\in I} \tau(d_{i_1}(1) d_{i_{\pi(1)}}(2) \ldots 
d_{i_q} (p-1) d_{i_{\pi(q)}} (p))\right| \le S_1S_2\ldots S_p\leqno (3.8)$$
where $S_k = S((d_i(k))_{i\in I}, p)$. More generally, for any $t\ge 1$, we have
$$\left\|\sum_{i_1\ldots i_q} d_{i_1}(1) d_{i_{\pi(1)}}(2) \ldots d_{i_q} (p-1) 
d_{i_{\pi(q)}} (p)\right\|_t \le \prod^p_{k=1} S((d_i(k))_{i\in I}, pt).\leqno 
(3.9)$$

\pf Note that (3.9) when $t=1$ obviously implies (3.8). We will prove (3.9) (for 
all $t\ge 1$) by induction on $q$. The case $q=1$ is very easy since it is well 
known that for all $t\ge 1$ 
$$\left\|\sum d_i(1) d_i(2)\right\|_t \le \left\|\left( \sum d_i(1) 
d_i(1)^*\right)^{1/2}\right\|_{2t} \cdot \left\|\left( \sum d_i(2)^* d_i(2) 
\right)^{1/2}\right\|_{2t}.\leqno (3.10)$$
Assume that (3.9) has been proved (for all $t\ge 1$) for the value $q-1$. Let us 
show that it also holds for $q$. By definition of $S^{nc}_q$, the partition of 
$[1,\ldots, 2q]$ into pairs associated to $\pi$ is non-crossing. This implies 
that this partition admits an interval $[k,k+1]$ as one of its blocks.
Moreover  if 
we delete this block the resulting partition of the remaining set (with the 
induced ordering) is still non-crossing. Let $x = \sum\limits_{i_1\ldots i_q} 
d_{i_1}(1) d_{i_{\pi(1)}}(2)\ldots d_{i_q}(p-1)d_{i_{\pi(q)}}(p)$. Thus we can 
write
$$x = \sum_\alpha a_\alpha \sum_{i\in I} d_i(k) d_i(k+1)b_\alpha$$
hence
$$\|x\|_t \le \left\|\sum_{i\in I} d_i(k) d_i(k+1)\right\|_{qt} \cdot C\leqno
(3.11)$$ where
$$C = \sup\left\{\left\|\sum a_\alpha Tb_\alpha\right\|_t ~\Big|~ \|T\|_{qt} \le 
1 \right\}.$$
Now, by the induction hypothesis we know that for any $s\ge 1$,  for
any 
$u$ with $\|u\|_\infty \le 1$,  we have
$$\left\|\sum_\alpha a_\alpha ub_\alpha\right\|_s \le C'
,$$
with $C' =\prod_{\xi\notin [k,k+1]} S((d_i(\xi))_{i\in I}, (p-2)s)$.

\n Thus the linear mapping $v$ defined by
$$v(y) = \sum_\alpha a_\alpha yb_\alpha$$
is bounded from $L_\infty$ into $L_s$ with norm $\le C'$. Since the partition 
corresponding to $\sum b_\alpha u a_\alpha$ obviously also is non-crossing, we 
also have the same bound for ${}^t{v(y)} = \sum b_\alpha ya_\alpha$, or 
equivalently we know that $v$ is bounded with norm $\le C'$ from $L_{s'}$ to 
$L_1$. By interpolation, for any $0<\theta <1$, it follows that $v$ is also 
bounded from $L_a$ to $L_b$ where
$${1\over a} = {1-\theta\over \infty} + {\theta\over s'},\quad {1\over b} = 
{1-\theta\over s} + {\theta\over 1}.$$
If we choose $s$ so that ${1\over s} = {1\over t} \big[1 - {1\over q}\big]$. 
Then imposing $b=t$, we find $\theta$ determined by $\theta(1-1/s) = {1\over b} 
- {1\over s} = {1\over t}-{1\over s}$. Then the value  of $a$ is given by 
${1\over a} = {\theta\over s'} = \theta \big(1- {1\over s}\big) = {1\over t} 
-  {1\over s} = {1\over qt}$. Thus we conclude that $v$ is bounded from
$L_{qt}$ to 
$L_t$ with norm $\le C'$. In other words, we have established that
$$C\le C'.$$
Note that $(p-2)s = 2(q-1)s = 2qt = pt$. Moreover, by (3.10) (applied in 
$L_{qt}$ instead of $L_t$) we have
$$\left\|\sum d_i(k)d_i(k+1)\right\|_{qt} \le C''
$$
with $C''=\left\|\left(\sum d_i(k) 
d_i(k)^*\right)^{1/2}\right\|_{pt}\cdot \left\|\left( \sum d_i(k+1)^*
d_i(k+1) 
\right)^{1/2}\right\|_{pt}$.

\n Hence we can finally deduce from (3.11) that $\|x\|_t \le CC'' \le C'C''$
and  since $(p-2)s = pt$ we find that $C'C''$ is less or equal to the right
side of  (3.9).\qed

\n{\bf Remark.} The analogs of Proposition 1.1 and Theorem 1.2 for the lattice
of non crossing partitions are proved in [Sp2]. Thus we can combine this
with the same argument as above if the function $\sigma\to \Phi(\sigma)$ is
supported by the set of non crossing partitions, and the resulting
constants will remain bounded when $p$ tends to $\infty$. However, we could
not find a  significant application of this idea.

\ve

\n {\bf \S 4. Applications to Harmonic Analysis}

 The results of this section can be viewed as a continuation
of a series of investigations devoted to Fourier series with coefficients
in a non-commutative $L_p$-space, such as \ e.g. [TJ, BP, LP, LPP, X].

As explained in the introduction, our main inequality applies to $p$-dissociate 
partitions $\Lambda = \bigcup\limits_{i\in I} \Lambda_i$ of a subset $\Lambda$ 
in a discrete group $G$.
The   inequality in Theorem 4.1 below is closely related (and partly
motivated) by the recent  papers [H1-2] on the so-called
$\Lambda(p)_{cb}$-sets, which are a certain  non-commutative version of
Rudin's classical
$\Lambda(p)$-sets (\cf [Ru]). The basic examples of such sets are the
$p$-dissociate ones.
However,
in the quest for the ``largest possible'' examples of sets  satisfying
such inequalities,  the next result turns out to be more  efficient
and more flexible (in particular in the analysis of
$\Lambda(p)_{cb}$-sets
constructed as random subsets of a given set), even
though its assumptions become   more complicated than the condition of
being $p$-dissociate.

\proclaim Theorem 4.1. Let $1 = \sum\limits_{j\in J} P_j$ be an orthogonal 
decomposition of the identity of $L_2(\tau)$. Let $p=2q$ be an even integer 
$>2$. Let $d = (d_i)_{i\in I}$ be a finite family in $L_p(\tau)$. We set 
$x^\omega = x^*$ if $q$ is odd and $x^\omega=x$ if $q$ is even.
Let $F$ be the set of all {\it injective\/} functions $g\colon \ [1,2,\ldots, 
q]\to I$. For any $g$ in $F$, we denote $x_g = d^*_{g(1)} d_{g(2)} d^*_{g(3)} 
\ldots d^\omega_{g(q)}$. We then define
$$N(d) = \sup_{j\in J} \hbox{ card}\{g\in F\mid P_jx_g\ne 0\}.$$
We then have
$$\left\|\sum_{i\in I} d_i\right\|_p \le \left[ (4N(d))^{1/p} + p \cdot {9\pi 
\over 8}\right] S(d,p).$$

\pf Since the argument is essentially the same as in [H2] modulo the 
combinatorics of \S 1, we will only sketch the proof.

\n Let $f = \sum\limits_{i\in I} d_i$. We have
$$\|f\|^q_p = \|f^*f\ldots f^\omega\|_2.$$
Developing this product as in \S 1 but with $n=q$ this time, $V = L_2(\tau)$ and 
$\varphi$ the product mapping, we obtain
$$f^*ff^*\ldots f^\omega = \Phi(\dot 0) - \sum_{\dot 0 < \pi \in P_q} \mu(\dot 
0,\pi)\Psi(\pi)\leqno (4.1)$$
where $\Phi(\sigma) = \sum\limits_{\pi(g)=\sigma} x_g$.

 \n  Using Sublemma 3.3 and 
(3.5) with $p$ replaced by $q$, we obtain (recall $\alpha= {3\pi\over 4}$)
$$\|\Psi(\pi)\|_2 \le \|f\|^{r_1(\pi)}_p (\alpha\|S\|_p)^{q-r_1(\pi)}. \leqno 
(4.2)$$
On the other hand, we can write 
$$\eqalign{\|\Phi(\dot 0)\|^2_2 &= \sum_{j\in J} \|P_j \Phi(\dot 0)\|^2_2\cr
&= \sum_{j\in J}\left\| \sum_{g\in F} P_jx_g\right\|^2_2\cr
&\le \sum_{j\in J} \sum_{g\in F} \|P_jx_g\|^2_2 N(d)\cr
&\le \sum_{g\in F} \|x_g\|^2_2 N(d)\cr
&\le \sum_{g_1,\ldots, g_q\in I} \tau(d^*_{g(1)}\ldots d^\omega_{g(q)} 
(d^\omega_{g(q)})^* \ldots d_{g(1)})\cdot N(d)}$$
hence by a special case of Lemma 3.6, we have
$$\|\Phi(\dot 0)\|^2_2 \le N(d) S(d,p)^{2q}.\leqno (4.3)$$
Thus, arguing as in the proof of Theorem 2.1,  we obtain finally combining
(4.1), (4.2) and (4.3):
$$\|f\|^q_p \le N(d)^{1/2} S(d,p)^q  + \sum_{0\le s<q} {q\choose s}(q-s)!
\|f\|^s_p  (\alpha S(d,p))^{q-s}$$
hence, if we now set $y = {S(d,p)\over \|f\|_p}$, we have
$$1\le N(d)^{1/2} y^q + \sum_{0\le s<q} {q\choose s}(q-s)! (\alpha y)^{q-s}.$$
We claim that $y \ge \min\big\{\big({1\over 2N(d)^{1/2}}\big)^{1/q}, {1\over 
3q\alpha}\big\}$. Indeed, if $y < (2N(d)^{1/2})^{-1/q}$ then, as in the proof
of Sublemma 2.3, we have
$$3/2 \le \int^\infty_0 (1+\alpha ty)^q e^{-t} \ dt$$
which yields, if $qy\alpha < 1$, that $3/2 \le ( 1-qy\alpha)^{-1}$ whence $y 
\ge {1\over 3q\alpha}$; otherwise $qy \alpha\ge 1$ which also implies $y \ge 
{1\over 3q\alpha}$.
Thus we conclude as announced that a fortiori we have
$$1/y \le (2N(d)^{1/2})^{1/q} + 3q\alpha.\eqno $$\qed

\proclaim Corollary 4.2. ([H2])~~Let $\Lambda\subset G$ be a subset of a
discrete  group $G$. Let $p=2q$ be an even integer $>2$. For any $t$ in $G$,
let 
$N_q(t,\Lambda)$ be the number of $q$-tuples $(t_1,\ldots, t_q)$ of mutually 
distinct elements of $\Lambda$ such that
$$t = t^{-1}_1 t_2t^{-1}_3t_4\ldots t^\omega_q.$$
We assume that
$$N_q(\Lambda) = \sup_{t\in G} N_q(t,\Lambda) < \infty.$$
Then, for any finitely supported family $a = (a_t)_{t\in\Lambda}$ in a 
non-commutative $L_p$-space associated to a semi-finite trace $T$, we have
$$\left\|\sum_{t\in\Lambda} \lambda(t) \otimes a_t\right\|_{L_p(\tau_G \times 
T)} \le \left[(4N_q(\Lambda))^{1/p} + p {9\pi\over 8}\right] S(a,p).$$

\pf We apply the previous result to $\tau = \tau_G \times  T$ so that 
$L_2(\tau) = L_2(\tau_G) \otimes_2 L_2(T)$ and to the $\bot$ decomposition 
$L_2(\tau) = \bigoplus\limits_{t\in G} H_t$ with $H_t = \lambda(t) \otimes 
L_2(T)$ and $I=\Lambda$. Clearly, if we set $d_t = \lambda(t) \otimes a_t$, 
$t\in\Lambda$ we find $N(d) \le N_q(\Lambda)$ and the result follows since 
$S(d,p) = S(a,p)$.

\qed \ve

\n {\bf \S 5. Tensor Products of Banach Spaces}

The main idea exploited above can also be used in a very abstract
setting, which we briefly indicate in this section.  Let 
$E_1,\ldots, E_p$ be Banach spaces and let $E_1\widehat\otimes \cdots 
\widehat\otimes E_p$ be their projective tensor product equipped with its 
projective norm denoted by $\|~~~\|_\wedge$ (see e.g.\ [DF]).

For each $k=1,2,\ldots, p$ consider a finite sum
$$f_k = \sum_{i\in I} d_i(k)$$
where $d_i(k)$ are elements of $E_k$.

Now let $(\vp_i)_{i\in I}$ be a sequence of independent $\pm 1$-valued random 
variables on a probability space $(\Omega, P)$ with $P(\vp_i = \pm 1) = 1/2$, as 
usual.\medskip

\n {\bf Note:}~~The family $(\vp_i)_{i\in I}$ is the Abelian counterpart of the 
family $(\lambda(g_i))_{i\in I}$ used above.

We wish to develop the tensor product 
$$f_1 \otimes \cdots \otimes f_p$$
in the Banach space $E_1\widehat\otimes\cdots \widehat\otimes E_p$. We will use 
the notation in \S 1 applied to the canonical multilinear mapping $\varphi\colon 
\ E_1\times\cdots \times E_p\to E_1 \widehat\otimes \cdots \widehat\otimes E_p$. 
Hence we have now
$$\Phi(\dot 0) = \sum_g d_{g(1)}(1) \otimes\cdots \otimes d_{g(p)}(p)$$
where the sum runs over all injective maps $g\colon \ [1,2,\ldots, p]\to I$. Let 
$\pi$ be a partition of $[1,\ldots, p]$. Using the random variables 
$(\vp_i)_{i\in I}$ instead of $(\lambda(g_i))_{i\in I}$ in the preceding 
section, it is easy to adapt the proof of Sublemma~3.2 to obtain the
following  result:

Let $A \subset [1,\ldots, p]$ be the union of the singletons of the partition 
$\pi$ (note that the cardinality of $A$ is at most $p-2$, unless $\pi=\dot
0$) and, as before, we let
$$\Psi(\pi) = \sum_{g\colon \ \pi(g)\ge \pi} d_{g(1)}(1) \otimes\cdots \otimes 
d_{g(p)}(p).$$
Then we have
$$f_1 \otimes \cdots \otimes f_p=\Phi(\dot 0) - \sum_{0<\pi} \Psi(\pi) \mu(\dot 0, 
\pi)
 $$
 and 
$$\|\Psi(\pi)\|_\wedge \le \prod_{k\in A} \|f_k\| \cdot \prod_{k\notin A} S_k$$
where
$$S_k = \left(\EE\left\|\sum_{i\in I} \vp_id_i(k)\right\|^p \right)^{1/p}.$$
We can now state the main result of this section.

\proclaim Theorem 5.1. With the above notation, we have
$$\left\|f_1\otimes\cdots\otimes f_p - \sum_{\sst g\colon \ [1,\ldots, 
p]\to I\atop \sst g~{\rm injective}} d_{g(1)} \otimes\cdots \otimes 
d_{g(p)}\right\|_\wedge
\le \sum_{\sst A\subset [1,\ldots,p]\atop \sst |A|\le p-2} \prod_{k\in A} 
\|f_k\| \cdot \prod_{k\notin A} S_k\cdot (p-|A|)!$$
In the particular case $E_1 = E_2 =\cdots= E_p=E$ we obtain:

\proclaim Corollary 5.2. Let $f  = \sum\limits_{i\in I} d_i$ be a finite sum in 
a Banach space $E$. Let $f^{\otimes p} = f \otimes\cdots \otimes f$ ($p$-times). 
Then
$$\left\|f^{\otimes p} - \sum_{\sst g\colon \ [1,\ldots, p]\to I\atop\sst g~{\rm 
injective}} d_{g(1)} \otimes\cdots\otimes d_{g(p)}\right\|_\wedge \le 
\sum_{0\le s\le p-2} {p\choose s} (p-s)! \|f\|^s S^{p-s}$$
where $S = \left(\EE\left\|\sum\limits_{i\in I} \vp_id_i\right\|^p 
\right)^{1/p}$.\ve

\centerline{\bf References}

\item{[A]} G. Andrews. The theory of partitions.
Cambridge Univ. Press, 1984.

\item{[B1]} D.\ Burkholder, Distribution function
 inequalities for martingales, Ann.
Probability 1 (1973) 19-42.

\item{[B2]} D.\ Burkholder, A geometrical characterization 
of Banach spaces
in which martingale
difference sequences are
 unconditional, Ann.\ Probab. 9 (1981), 997-1011.

\item{[B3]} D. Burkholder. Boundary value problems and sharp inequalities
for martingale transforms. Ann. Prob. 12 (1984) 647-702.

\item{[B4]} D. Burkholder. Explorations in martingale theory
and its applications. (Ecole d'\'et\'e  de Probabilit\'es de Saint-Flour 
XIX, 1989) Springer Lecture Notes 1464 (1991) 1-66.

\item{[BG]} D.\ Burkholder and  R. F.
   Gundy. Extrapolation and interpolation of
quasi-linear operators on martingales. Acta Math. 124 (1970) 249-304.

\item{[BP]} O. Blasco and A. Pe{\l}czy{\'n}ski. Theorems of Hardy and Paley
for vector-valued analytic functions and related classes of Banach spaces. 
Trans. Amer. Math. Soc.
323 (1991)  335-367. 

\item{[BS]} P. Biane and  R. Speicher. Stochastics calculus with respect to
free
Brownian motion and analysis on Wigner space.  Probab. Theory 
Related Fields 112
(1998)   373-409.

\item{[BSp]} M. Bo\.{z}ejko and  R.\ Speicher. An example of
a generalized brownian motion. Comm. Math.
Physics. 137 (1991) 519-531.

\item{[Bu1]} A.\ Buchholz. Norm of convolution by operator-valued functions on 
free groups.  Proc. Amer. Math. Soc., to appear.

\item{[Bu2]} A.\ Buchholz. Operator Khintchine inequality for $q$-Gaussian 
random variables,  preprint 1998, to appear in Math. Ann.

\item{[DF]}  A. Defant, K. Floret. Tensor norms and operator 
ideals.
North-Holland, Amsterdam, 1993.

\item{[HT]} U. Haagerup and S. Thorbj\o rnsen. Random matrices and $K$-theory
for
exact $C^*$-algebras. Odense University preprint, 1998.

\item{[H1]} A. Harcharras. Analyse de Fourier, multiplicateurs de Schur
sur $S_p$ et ensembles $\Lambda(p)_{cb}$ non commutatifs.
Comptes Rendus Acad. Sci. Paris 326
(1998) 845-850.

\item{[H2]} A. Harcharras. Fourier analysis, Schur multipliers on $S_p$
and non-commutative $\Lambda(p)$-sets. Studia Math. To appear.

\item{[K]} G.\ Kreweras.  Sur les partitions non-crois\'ees d'un cycle. Discrete 
Math. {\bf 1} (1972), 333--350.

\item{[LP]} F.\ Lust-Piquard, In\'egalit\'es de Khintchine dans $C_p$
$(1<p<\infty)$, C.R.\ Acad.\ Sci. Paris 303 (1986), 289--292.

\item{[LPP]} F.\ Lust-Piquard, G.\ Pisier, Noncommutative Khintchine and
Paley inequalities, Arkiv
f\"or Mat. 29 (1991), 241--260.

\item{[LS]}   M.  Longnecker and  R. J. Serfling. 
        Moment inequalities for $S_n$ under general dependence restrictions,
        with applications.
        Z. Wahrsch. verw. Geb.  43 (1978) 1-21.

\item{[Pa]} R. E. A. C. Paley, A remarkable series of orthogonal functions (I),
Proc. London Math. Soc. 34 (1932) 241-264.

\item{[PX]} G. Pisier and Q. Xu.
Non-commutative martingale inequalities.
Comm. Math. Physics 189 (1997) 667-698.

\item{[R1]} G.C.\ Rota. On the foundations of combinatorial theory I:\ Theory of 
M\"obius functions. Z.\ Warschein.\ Verw.\ Geb. {\bf 2} (1964), 340--368.

\item{[R2]} G.C.\ Rota. Th\'eorie combinatoire des invariants classiques. 
Seminar Notes 76/77. Strasbourg (France).

\item{[Ru]} W. Rudin. Trigonometric series with gaps. J. Math. Mech. 9 (1960) 203-228.

\item{[Se]} R. J. Serfling. Probability inequalities and convergence
properties
for sums of multiplicative random variables. Unpublished preprint.
(FSU Statistics report M151, February 1969).

\item{[Sp1]} R.\ Speicher. A new example of ``independence and white noise".
Probab. Th. rel. Fields. 84 (1990) 141-159.

\item{[Sp2]} R.\ Speicher. Multiplicative functions on the lattice of 
non-crossing 
partitions and free convolution. Math.\  Ann. {\bf 298} (1994), 611--628.

\item{[St]} W. Stout. Almost sure convergence. Academic Press, New-York, 1974.

\item{[TJ]} N. Tomczak-Jaegermann. The moduli of 
convexity and smoothness and the Rademacher averages of trace class
$S_p$. Studia Math. 50 (1974) 163-182.

\item{[VDN]} D.V. Voiculescu, K.J. Dykema, A. Nica, Free Random variables,
CRM Monograph Series,
Vol.1, Centre de Recherches Math\'ematiques, Universit\'e Montr\'eal,

\item{[X]} Q. Xu. Analytic functions with values in lattices and
symmetric spaces of measurable operators. 
Math. Proc. Cambridge Philos. Soc. 109 (1991),
no. 3, 541-563.

\bye